\documentclass[12pt]{article}

\usepackage{sectsty}
\sectionfont{\color{dark}} 
\subsectionfont{\color{dark}} 
\subsubsectionfont{\color{light}} 
\usepackage{cancel,enumerate}
\usepackage[affil-it]{authblk}

\usepackage[intlimits]{amsmath}
\usepackage{amstext, amssymb}
\usepackage[mathscr]{euscript}
\usepackage{eufrak}
\usepackage[mathscr]{euscript}

\usepackage{fancyvrb}

\DefineVerbatimEnvironment%
{MapleInputs}{Verbatim}{formatcom=\color{MapleRed}}

\newcommand{\st}{\,|\,}

\usepackage[strict]{changepage}
\usepackage{dsfont}
\usepackage{tikz}
\usetikzlibrary{intersections}
\usetikzlibrary{snakes}
\usetikzlibrary{arrows}
\usetikzlibrary{shapes}
\usetikzlibrary{backgrounds}
\usetikzlibrary{shadows}
\usetikzlibrary{decorations.pathreplacing}
\usetikzlibrary{decorations.markings}
\usepackage{environ}
\NewEnviron{sidenote}{
\begin{center}

  \noindent  
  \begin{tikzpicture}[font=\small]
      \node [text width=0.9\linewidth,align=justify,draw,drop shadow,fill=white,inner sep=10pt]{%
           \BODY
        };
   \end{tikzpicture}
\end{center}
}

\usepackage {
	amsmath,
	amssymb,
	amsthm,
	amsxtra,
	appendix,
	epsfig,
	fourier-orns,
	graphics,
	longtable,
	needspace,
	makeidx,
	mdwlist,
	multicol,
	psfrag,
	pstricks,
	rotating,
	setspace,
	sidecap,
	stmaryrd,
	subfloat,
	tabularx,
	textcomp,
	tocloft,
	url,
	wrapfig
}

\usepackage[mathscr]{euscript}
\usepackage{etoolbox}
\usepackage[linesnumbered,boxed,vlined]{algorithm2e}
\SetAlCapSkip{2ex}
\SetAlgoVlined
\SetAlgoInsideSkip{medskip}

\usepackage{tikz}
\usetikzlibrary{positioning,calc,decorations.pathreplacing}
\usepackage[font={small}]{caption}	
\usepackage{lipsum}
\usepackage{subcaption}
\usepackage{etex}
\usepackage{mathtools}
\usepackage{needspace}

\newif\ifisappendix

\allowdisplaybreaks

\AtBeginDocument{%
    \mathchardef\mathcomma\mathcode`\,
    \mathcode`\,="8000 
    \mathchardef\mathsemicolon\mathcode`\;
    \mathcode`\;="8000 
}
{\catcode`,=\active
  \gdef,{\mathcomma\discretionary{}{}{}}
\catcode`;=\active
  \gdef;{\mathsemicolon\discretionary{}{}{}}
}

\newtheoremstyle{coloured}{}{}{}{}{\color{light}\bfseries}{.}{5pt plus 1pt minus 1pt}{}
\theoremstyle{coloured}

\newtheorem{definition}{Definition}[section]
\newtheorem{example}{Example}[section]

\newcommand{\cbrac}[1]{\{ #1 \}}

\newcommand{\HH}{{\mbox{\footnotesize\sf H}}}
\newcommand{\TT}{{\mbox{\footnotesize\sf T}}}
\newcommand{\RR}{{\mbox{\footnotesize\sf R}}}

\newcommand{\ii}{\mathbf{i}}
\newcommand{\jj}{\mathbf{j}}

\renewcommand{\vec}{\mathbf}

\newcommand{\1}{\mathbf{1}}
\newcommand{\2}{\mathbf{2}}

\newcommand{\cS}{\mathcal{S}}

\makeatletter
\newcommand{\superimpose}[2]{%
  {\ooalign{$#1\@firstoftwo#2$\cr\hfil$#1\@secondoftwo#2$\hfil\cr}}}
\makeatother

\newcommand{\cW}{\mathcal{W}}
\newcommand{\dee}{{\rm d}}
\newcommand{\mm}{\vec{m}}

\definecolor{dark}{RGB}{0, 133, 202 }
\definecolor{light}{RGB}{112, 155, 230}

\renewcommand{\emptyset}{\varnothing}

\begin{document}

\title{\color{dark} Playing Several Patterns Against One Another}

\author[1]{Jan Vrbik}
\author[2]{Paul Vrbik}
\affil[1]{Department of Mathematics and Statistics\\ Brock University, Canada}
\affil[2]{School of Mathematical and Physical Sciences, University of Newcastle, Australia}

\date{\today}
\maketitle

\begin{abstract}
We revisit the game in which each of several players chooses a pattern and then a coin is flipped repeatedly until one of these patterns is generated. In particular, we demonstrate how to compute the probability of any one player winning this game, and find the distribution of the game's duration. Our presentation is an extension (and perhaps a simplification) of the results of Blom and Thornburn \cite{ref}. 
\end{abstract}

\section{Introduction}

Blom and Thornburn in \cite{ref} outline an experiment where a coin or a die is flipped or rolled until a particular pattern (such as \TT\TT\HH~or 123) is generated. Different patterns can then be `played against' one another by designating the winning pattern as the first one that is generated. In \cite{ref} however, it is assumed that the coin is fair and that the patterns are all of the same length. In this paper we remove these restrictions; we also aim our presentation at a wider audience by simplifying the key definitions and ideas.\medskip

We build our patterns out of only two possible characters (\TT and \HH) because the extension to using three or more is trivial.

\section{Pattern Generation}\label{Fri24April527}

Consider performing a sequence of independent \textsc{trials} by flipping a potentially biased coin. Each trial results in a `head' (denoted $\HH$ in this article) with a probability $p$, or `tail' ($\TT$) with probability $q\equiv 1-p$. Flips continue until a specific \textsc{pattern} of $k$ consecutive outcomes (e.g. $\HH\HH\TT\HH\TT$) is generated. In this article, we use calligraphic capital letters (such as $\cS$) to denote a pattern and small letters (such as $s_{1}s_{2}\cdots s_{k}$) for its individual symbols.

\begin{definition}[First-Time Generation]\label{Tue14Apr957}
Let $f_{i}$ denote the probability that completing a specific pattern $\cS$ of length $k$ \emph{for the first time} will happen at Trial $i$. Thereby $f_{i}=0$ when $i<k$, including $i=0$. The probability generating function (PGF) of the corresponding $f_i$ sequence is 
$$
F(z)\equiv \sum_{i=0}^{\infty }f_{i}\cdot z^{i}.
$$
\end{definition}

The expected (or mean) value of the number of trials to generate the pattern is then given by $F^{\prime}(z=1)$. 

\subsection{A key formula}

Assume now that trials are repeated indefinitely and generate an ever increasing number of occurrences of pattern $\cS$.

\begin{definition}\label{T14A1006}
Let $u_{i}$ denote the probability that Pattern $\cS$ is completed at Trial $i$ but not necessarily for the first time. Here $u_{i}=0$ when $0<i<k$, but $u_{0}=1$ (we give the rationale below). Let $U(z)$ be the generating function of the $u_{0},\, u_{1},\, u_{2},\, \ldots$ sequence. Note that this generating function is \emph{not} a PGF since the $u_{i}$ probabilities do not add up to $1$.
\end{definition}

Definition \ref{T14A1006} requires an important \emph{proviso}: consecutive occurrences of the pattern are  \emph{not} allowed to overlap. Each time Pattern $\cS$ is generated, the process is reset and the next completion has to start `from scratch', with none of the symbols of any one occurrence being available to build the next one. This implies that not every subsequence of the $s_1s_2 \cdots s_k$ symbols (we call it String $\cS$, to emphasize the difference) counts as a completion of Pattern $S$. \medskip

Let us now assume that exactly $n$ trials have been completed (this becomes our \textsc{sample space}) and let us expand $u_{n}$ according to the trial of the first occurrence of Pattern $\cS$ (thus defining a \textsc{partition} of the sample space). Using the total-probability formula, we deduce
\begin{equation}
u_{n}=f_{0}u_{n}+f_{1}u_{n-1}+\cdots+f_{n}u_{0}  \label{uf}
\end{equation}
where each term of the right-hand side, say $f_{i}u_{n-i}$, represents the probability of completing the \emph{first} occurrence of $\cS$ at Trial $i$, followed by generating yet another (not necessarily second) occurrence of $\cS$ at the end of the remaining $n-i$ trials. The first term with $i=0$ is always zero; the last one accounts for the possibility of \emph{no prior} occurrence of $\cS$, thus explaining the necessity of choosing $u_{0}=1$. Since the terms represent all possibilities of what can happen to complete $\cS$ at Trial $n$, their sum must equal the left-hand side.

Notice that the RHS of (\ref{uf}) is equal to the coefficient of $z^{n}$ in the expansion of $F(z)\cdot U(z)$, known as the \textsc{convolution} of the two sequences. It is also important to realize that (\ref{uf}) is correct for $n\geq 1$ but not when $n=0$. Multiplying each side of (\ref{uf}) by $z^{n}$ and summing from $1$ to $\infty$ thus yields
$$
U(z)-1=F(z)\cdot U(z)
$$
where the `$-1$' is to account for the missing $u_{0}$ on the LHS. Solving for $F(z)$ yields
\begin{equation}
F(z)=\frac{U(z)-1}{1+(U(z)-1)}=\left(1+\frac{1}{U(z)-1}\right)^{-1} \label{fu}
\end{equation}

\subsection{Pattern Overlaps}

For our presentation to be concise and readable, we introduce some operations on strings. Namely, we frequently need to know if and when the last few symbols of a pattern coincide with the first few symbols of another pattern. 

\begin{example}\label{Sat4Jul1139}
$\TT\TT\TT\HH\TT\TT\TT$ overlaps with $\TT\TT\HH\TT\TT\TT\TT\HH\TT$ \emph{only} at shift 1, 2 and 6 because (red indicates the overlap)
\begin{align*}
\mbox{Shift 1} && \TT\,\TT\,\TT\,\HH\,\TT\,\TT\, &{\color{red}\TT}  \\[-.5em]
&& &{\color{red}\TT}\,\TT\,\HH\,\TT\,\TT\,\TT\,\TT\,\HH\,\TT \\
\mbox{Shift 2} && \TT\,\TT\,\TT\,\HH\,\TT\, &{\color{red}\TT\, \TT}  \\[-.5em]
&& &{\color{red}\TT\,\TT\,}\HH\,\TT\,\TT\,\TT\,\TT\,\HH\,\TT \\
\mbox{Shift 6} && \TT\,&{\color{red}\TT\,\TT\,\HH\,\TT\,\TT\,\TT}  \\[-.5em]
&& &{\color{red}\TT\,\TT\,\HH\,\TT\,\TT\,\TT}\,\TT\,\HH\,\TT 
\end{align*}
are the only three possibilities. 

And similarly, reversing the order now, $\TT\TT\HH\TT\TT\TT\TT\HH\TT$ overlaps with $\TT\TT\TT\HH\TT\TT\TT$ overlap \emph{only} at shift 1 and 5:
\begin{align*}
\mbox{Shift 1} && \TT\,\TT\,\HH\,\TT\,\TT\,\TT\,\TT\,\HH\, &{\color{red}\TT}  \\[-.5em]
&& &{\color{red}\TT}\,\TT\,\TT\,\HH\,\TT\,\TT\,\TT \\
\mbox{Shift 5} && \TT\,\TT\,\HH\,\TT\, &{\color{red}\TT\,\TT\,\TT\,\HH\,\TT}  \\[-.5em]
&& &{\color{red}\TT\,\TT\,\TT\,\HH\,\TT}\,\TT\,\TT \\
\end{align*}
\end{example}

\begin{definition}[$\boxdot$]
Let $\boxdot$ be a non-commutative function with two \emph{string} arguments which returns the set of shifts at which the first patten overlaps with the second one. Namely, when $\cS = s_1 \cdots s_k$ and $\cW = w_1 \cdots w_m$ then
$$\cS \boxdot \cW := 
\cbrac{ i \st s_{k-i+1} \cdots s_k = w_1 \cdots w_i }.
$$
This set is empty when there are no such overlaps. Note that this definition implies $i \leq \max(k,m)$.
\end{definition}

\begin{definition}[$\circledcirc$]
Let $\circledcirc$ be a non-commutative function with two \emph{string} arguments which returns the largest shift possible for an $\cS$ and $\cW$ overlap.
$$\cS \circledcirc \cW := \max( i \in \cS \boxdot \cW  )$$
where $\cS \circledcirc \cW := 0$ when $\cS \boxdot \cW = \emptyset$.
\end{definition}

\begin{definition}[$\sqcap$]
Let $\sqcap$ be a non-commutative function with two \emph{string} arguments which returns the substring corresponding to the longest $\cS$ and $\cW$ overlap:
$$\cS \sqcap \cW := w_{1} \cdots w_{\cS \circledcirc \cW}$$
($\cS \sqcap \cW := s_{k+1-\cS \circledcirc \cW} \cdots s_k$ is an equivalent definition). Note that $\cS \sqcap \cW$ is a zero-length string when $\cS \circledcirc \cW = 0$.
\end{definition}


\begin{example}\label{FirstExample}
Assuming that 
$\cS=\TT\TT\TT\HH\TT\TT\TT$ and 
$\cW= \TT\TT\HH\TT\TT\TT\TT\HH\TT$ (the strings of Example \ref{Sat4Jul1139}) we get
\begin{align*}
\cS\boxdot \cW &= \{1,2,6\} 
&&& \cW\boxdot \cS &= \{1,5\} \\
\cS\circledcirc \cW &= 6 
&&& \cW\circledcirc \cS &= 5 \\
\cS\sqcap \cW &= \TT\TT\HH\TT\TT\TT 
&&& \cW\sqcap \cS &= \TT\TT\TT\HH\TT \\
\cS\boxdot \cS &= \{1,2,3,7\} 
&&& 
\cW\boxdot \cW &= \{1,4,9\}
\end{align*}
\end{example}

\subsection{Finding $U(z)$ and $F(z)$}

Using the same sample space of $n$ trials, consider the probability of obtaining the \emph{symbols} $s_{1}\cdots s_{k}$ at the last $k$ of these $n$ trials. This does \emph{not} necessarily mean that Pattern $\cS$ has been completed at Trial $n$, since such an occurrence of \textsc{String} $\cS$ (as we call it, to differentiate from Pattern $\cS$) may overlap with an earlier completion of the actual pattern and therefore not count as another completion of Pattern $\cS$.\medskip

Denote the probability of generating String $\cS$ in (any) $k$ consecutive trials by $P_{\cS}$. It is clearly given by
\begin{equation}\label{Su19Apr1250}
P_{\cS} = p^{m}(1-p)^{k-m}
\end{equation}
where $m$ and $k-m$ are (resp.) the number of $\HH$ and $\TT$ symbols in $\cS$.

Now, we expand the probability of \emph{String} $\cS$ occurring at Trial $n$ (which is equal to $P_{\cS}$ for any $n\geq k$) according to the trial at which Pattern $\cS$ has been been completed during the last $k$ trials, thus: 
\begin{equation}
P_{\cS}=\sum_{\ell \epsilon \cS\boxdot \cS%
}u_{n-k+\ell }\cdot P_{s_{\ell+1}\cdots s_k}  \label{Pq}
\end{equation}
where $s_{k+1} \cdots s_k = \varepsilon$ (the empty string) and $P_\varepsilon = 1$. In this context, one must understand that when the last $k$ trials contain the symbols $s_{1}\cdots s_{k}$, Pattern $\cS$ must have been completed at \emph{exactly one} of these trials, and that only Trials $n-k+\ell $ with $\ell \in \cS \boxdot \cS$ are eligible.

Multiplying each side of (\ref{Pq}) by $z^{n}$ and summing over $n$ from $k$ to $\infty $ yields
\begin{align*}
\frac{P_{\cS}\cdot z^{k}}{1-z} 
&= \sum_{\ell \epsilon \cS \boxdot \cS}P_{s_{\ell + 1}\cdots s_k}\cdot z^{k-\ell }\sum_{n=k}^{\infty }u_{n-k+\ell }\cdot z^{n-k+\ell } \\
&= \sum_{\ell \epsilon \cS\boxdot \cS}P_{s_{\ell + 1}\cdots s_k}\cdot z^{k-\ell }\sum_{m=0}^{\infty }u_{m+\ell }\cdot z^{m+\ell } \\
&=\left( U(z)\overset{}{-}1\right) \cdot \sum_{\ell \epsilon \cS \boxdot \cS}P_{s_{\ell + 1}\cdots s_k}\cdot z^{k-\ell }
\end{align*}
since $1\leq \ell \leq k$ and $u_{1}=\cdots=u_{k-1}=0$. Thus, in the $\sum_{m=0}^{\infty }u_{m+\ell }\cdot z^{m+\ell }$ summation, we are always missing $u_{0}$, but are including the rest of the $U(z)$ expansion.

Based on (\ref{fu}), we get
\begin{equation}\label{OneStar}
F(z)=\left( 1+(1-z)\cdot \frac{\sum_{\ell \epsilon \cS\boxdot\cS}P_{s_{\ell+1} \cdots s_k}\cdot z^{k-\ell }}{P_{\cS}\cdot z^{k}}\right) ^{-1}
\end{equation}
and can thereby find the PGF of the number of trials to generate, for the first time, any given pattern.

\begin{example}
Let $\cS = \TT\TT\HH\TT\TT\TT\TT\HH\TT$, as in the previous example. We get
$$
F(z)=\left( 1+(1-z)\cdot \frac{p^{2}q^{6}z^{8}+pq^{4}z^{5}+1}{p^{2}q^{7}z^{9}}\right) ^{-1}.
$$
By expanding this PGF in powers of $z$, we can extract the individual $f_{i}$ probabilities (as coefficients of $z^{i}$). For instance $f_{18}=p^{2}q^{7}(1-pq^{4}-p^{2}q^{6}-p^{2}q^{7})$ is the probability that $\cS$ is generated, for the first time, in exactly $18$ trials.
\end{example}

\section{Playing $m$ Patterns}

Let us now assume that we have a collection of $m$ patterns $\cS_1,\, \ldots,\, \cS_{m}$,
none of them being a substring of any other. The question is: What is the probability that $\cS_i$ (from this set) is completed before any of the other patterns, thereby `beating' them and `winning' this game?

\begin{definition}
Let $x_{\ii,n}$ be the probability of a \emph{first time} occurrence of $\cS_i$ (or `Pattern $\ii$' for short) at Trial $n$, without being preceded by any of the competing patterns. In other words, this is the probability of $\cS_i$ beating all other patterns in exactly $n$ trials. For consistency we set $x_{\ii,0} = 0$.
\end{definition}

\subsection{Head-start probabilities}

Recall that the string denoted $\cS_i$ consists of the same symbols as the pattern now denoted $\ii$. 
Completion of String $\cS_i$ does not necessarily imply completion of Pattern $\ii$.
\begin{definition}
Let $f_{\ii,n}$ denote the probability of a first-time completion of Pattern $\ii$ at Trial $n$ (what was $f_n$ of Section \ref{Fri24April527}).
Furthermore, let $f_{\ii,n}^{(\jj)}$ be the same probability, this time assuming that Pattern $\jj$ has just been completed (prior to starting our trials) and that we \emph{are allowed} to use any of its symbols to help us generate Pattern $\ii$.
\end{definition}

Obviously Pattern $\ii$ cannot utilize more than $\cS_j\circledcirc $$\cS_i$ symbols of Pattern $\jj$; having $\cS_j\sqcap $$\cS_i$ as our `head start' thus results in the same $f_{\ii,n}^{(\jj)}$ probabilities. Let $F_{\ii|\jj}(z)$ denote their PGF.


To generate $\cS_i$ from scratch, one must first generate $\cS_j\sqcap $$\cS_i$ and then, independently (which implies convolution of the corresponding PGFs) the rest of the pattern. Denoting the PGF of $S_i \sqcap S_j$ by $F_{\jj\sqcap \ii}(z)$ we get
$$
F_{\ii}(z)=F_{\jj\sqcap \ii}(z)\cdot F_{\ii|\jj}(z),
$$
easily solved for $F_{\ii|\jj}(z)$ . Note that $F_{\ii}(z)$ and $F_{\jj \sqcap \ii}(z)$ can be found by utilizing Equation (\ref{OneStar}).

\subsection{Solving for $x_{\ii,n}$}

Assuming again that a fixed number of $n$ independent trials have been performed, we can expand $f_{\ii,n}$ according to the trial (say $\ell^{\rm th}$) in which Pattern $\jj$ has won the game. Notice that going over all values of $\ell $ and $\jj$ creates an obvious partition of our samples space. So,
\begin{equation}
f_{\ii,n}=\sum_{\jj\neq \ii}^{m}\sum_{\ell =1}^{n}x_{\jj,\ell }\cdot f_{\ii,n-\ell }^{(\jj)}+x_{\ii,n}
\label{key}
\end{equation}
as Pattern $\ii$ can only win at Trial $n$.

Multiplying (\ref{key}) by $z^{n}$ and summing over $n$ from $0$ to $\infty$ converts this infinite set of equations into a single statement involving the corresponding generating functions.  Namely:
\begin{align}
F_{\ii}(z) 
&=\sum_{\jj\neq \ii}^{m}X_{\jj}(z)\cdot F_{\ii|\jj}(z)+X_{\ii}(z) \notag\\ 
&=\sum_{\jj=1}^{m}X_{\jj}(z)\cdot F_{\ii|\jj}(z) \label{rec}
\end{align}
by setting $F_{\ii|\ii}(z)=1$. 

Dividing (\ref{rec}) by $F_{\ii}(z)$ converts it into\begin{equation}
1=\sum_{\jj=1}^{m}X_{\jj}(z)\cdot F_{\jj\sqcap \mathbf{i}}(z)^{-1}  \label{rec2}
\end{equation}
which has a simple solution, given by
\begin{equation}
\mathbf{X}=\mathbb{F}^{-1}\1  \label{fin}
\end{equation}%
where 
\begin{enumerate}[i.]
\item $\mathbf{X}$ is a vector of the $X_{\1}(z),\,\ldots,\,X_{\mathbf{m}}(z)$ generating functions,
\item $\1$ is a length-$m$ vector with each component equal to $1$, and 
\item $\mathbb{F}$ is the following matrix
$$
\left[\begin{array}{cccc}
\dfrac{1}{F_{\1 \sqcap \1}(z)} & \dfrac{1}{F_{\2 \sqcap \1}(z)} & \cdots & \dfrac{1}{F_{\mm \sqcap \1}(z)} \\[1em]
\dfrac{1}{F_{\1 \sqcap \2}(z)} & \dfrac{1}{F_{\2 \sqcap \2}(z)} & \cdots & \dfrac{1}{F_{\mm \sqcap \2}(z)} \\[1em]
\vdots & \vdots & \ddots & \vdots \\[1em]
\dfrac{1}{F_{\1 \sqcap \mm}(z)} & \dfrac{1}{F_{\2 \sqcap \mm}(z)} & \cdots & \dfrac{1}{F_{\mm \sqcap \mm}(z)}
\end{array}\right]
$$
That is $\mathbb{F}_{i,j} = {1}/{F_{\jj \sqcap \ii}}$ with $\ii$ specifying the \emph{row} index and $\jj$ the \emph{column} index (note the index reversal).
\end{enumerate}

To find the probability of Pattern $\ii$ winning the game \emph{in any number of trials}, we need to evaluate $X_{\ii}(z)$ at $z=1$. Since the substitution leads to an indefinite expression, we have to replace it by the corresponding $z\to 1$ limit (an easy task for a modern computer).

\begin{example}\label{We15Apr455}
For the patterns $\cS=\TT\TT\TT\HH\TT\TT\TT$ and $\cW= \TT\TT\HH\TT\TT\TT\TT\HH\TT$ of Example \ref{FirstExample} we have the following $\mathbb{F}$:
$$
\left[ 
\begin{array}{cc}
1+(1-z)\cdot \dfrac{pq^{5}z^{6}+pq^{4}z^{5}+pq^{3}z^{4}+1}{pq^{6}z^{7}} & 
1+(1-z)\cdot \dfrac{pq^{3}z^{4}+1}{pq^{4}z^{5}} \\[1em] 
1+(1-z)\cdot \dfrac{pq^{4}z^{5}+pq^{3}z^{4}+1}{pq^{5}z^{6}} & 1+(1-z)\cdot 
\dfrac{p^{2}q^{6}z^{8}+pq^{4}z^{5}+1}{p^{2}q^{7}z^{9}}%
\end{array}%
\right]
$$
resulting, with the help of \eqref{fin}, in 
\begin{equation}
\lim_{z\to 1}\mathbf{X}=\left[ 
\begin{array}{c}
\dfrac{1-p^{2}q^{3}}{1+p^{2}q-p^{3}q^{3}} \\[1.2em]
\dfrac{p^{2}q(1+q^{3})}{1+p^{2}q-p^{3}q^{3}}%
\end{array}
\right].\label{xxx}
\end{equation}
When $p=\frac{1}{2}$, this limit evaluates to $87.32\%$ and $12.68\%$ respectively, corresponding to the probability of $\cS$ and $\cW$ winning the game using a fair coin.
\end{example}
\subsection{The game's duration}

The PGF of the number of trials to complete the game is clearly given by
\begin{equation}
\sum_{\jj=1}^{m}X_{\jj}(z)  \label{dur}
\end{equation}
which can be readily expanded in $z$ to yield the corresponding probability of the game ending in exactly $i$ trials (the coefficient of $z^i$). To get the game's expected duration, we need to differentiate (\ref{dur}) with respect to $z$ and evaluate the result at $z=1$ (again, this will require taking the $z\to 1$ limit). Similarly, we can find the variance, skewness, and so on, of the corresponding distribution. 

\begin{example}
The game of Example \ref{We15Apr455} has expected duration given by
\begin{equation}
\lim_{z\to 1}\frac{\dee \left( X_{\1}(z)+X_{\2}(z)\right) }{\dee z}=\frac{1+q^{3}(1-q^{3})+p^{3}q^{6}+p^{2}(1+p)q^{9}}{pq^{6}(1+p^{2}q-p^{3}q^{3})}.  \label{mean}
\end{equation}%
This yields the average number of $128.3$ trials when $p=\frac{1}{2}$. Similarly, we get $\pm122.0$ trials for the standard deviation.
\end{example}

\subsection{Alternate solution}

There is an interesting way to bypass (\ref{fin}) when finding the probability of a Pattern winning the game.

Imagine playing the game repeatedly, without ever stopping; 
let $y_{\ii,n}$ denote the probability that Pattern $\ii$ wins the game at Trial $n$.
Provided we have played sufficiently many games to have reached the so-called \textsc{equilibrium}, the probability of Pattern $\ii$ winning a game is the same for all $n$, namely
$$y_{\ii} = \lim_{n \to \infty} y_{\ii,n}.$$ 
The probability of a game ending at Trial $n$ is also the same for all (sufficiently large) $n$, and is equal to $\sum_{\ii=1}^{m}y_{\ii}$. The expected duration of a single game must then be the reciprocal of this value, namely
\begin{equation}
\frac{1}{\sum_{\ii=1}^{m}y_{\ii}}. \label{Y2}
\end{equation}

Similarly, the probability of Pattern $\jj$ winning a game is
\begin{equation}
\frac{y_{\jj}}{\sum_{\ii=1}^{m}y_{\ii}}. \label{Y1}
\end{equation}
This will always agree with $\lim_{z \to 1} X_i(z)$ of Equation \eqref{fin}.

\subsubsection{Solving for $y_{\jj}$}

To find these $y_{\jj}$ probabilities, let $k$ be the length of Pattern $\jj$ and consider $k$ consecutive trials after reaching equilibrium; label these Trial 1 through Trial $k$. Then expand the probability of generating the \emph{symbols} of Pattern $\jj$ (i.e. of String $\cS_j$) in these $k$ trials, broken down according to the pattern which won the game at Trial $\ell$ (\emph{exactly one} pattern must have won during these $k$ trials).

Remembering that the probability of generating String $\cS_\jj = s_1 \cdots s_k$ is $P_{\cS_j}$, we get by summing over all patterns and all `eligible' trials:
\begin{equation}
P_{\cS_j}=\sum_{\ii=1}^{m}y_{\ii}\sum_{\ell\epsilon \cS_i\boxdot \cS_j}P_{s_{\ell+1} \cdots s_{k}}  \label{Pyi}
\end{equation}
where $P_{s_{\ell+1} \cdots s_{k}}$ is the probability of the corresponding string being completed during the last $k-\ell$ trials.

Doing this for each of the $m$ patterns provides $m$ linear equations for $y_{\1}$, $y_{\2},$ \ldots, $y_{\mathbf m}$. Solving them and substituting back into (\ref{Y1}) and (\ref{Y2}) yields, respectively, the probability of Patter $\jj$ winning, and the expected game's duration.

\begin{example}
For $\cS=\TT\TT\TT\HH\TT\TT\TT$ and $\cW= \TT\TT\HH\TT\TT\TT\TT\HH\TT$, the set of equations (\ref{Pyi}) reads
$$
\left[ 
\begin{array}{cc}
pq^{5}+pq^{4}+pq^{3}+1 & pq^{5}+q^{2} \\ 
p^{2}q^{6}+p^{2}q^{5}+pq^{2} & p^{2}q^{6}+pq^{4}+1
\end{array}
\right] \cdot \left[ 
\begin{array}{c}
y_{\1} \\ 
y_{\2}
\end{array}
\right] =\left[ 
\begin{array}{c}
pq^{6} \\ 
p^{2}q^{7}
\end{array}
\right] 
$$
resulting in 
\begin{align*}
y_{\1} 
&= \frac{pq^{6}(1-p^{2}q^{3})}{1+q^{3}(1-q^{3})+p^{3}q^{6}+p^{2}(1+p)q^{9}} \\
y_{\2} 
&= \frac{p^{3}q^{7}(1+q^{3})}{1+q^{3}(1-q^{3})+p^{3}q^{6}+p^{2}(1+p)q^{9}}
\end{align*}
which leads to (\ref{mean}) and (\ref{xxx}) via (\ref{Y2}) and (\ref{Y1}).
\end{example}

The advantage of this approach is that it bypasses the use of PGFs, allowing us to deal directly with numbers, rather than polynomials. This difference matters greatly when inverting matrices. Its shortcoming is its inability to compute more than the probabilities of winning the game and the expected value of the game's duration.

\subsection{Three Characters and Beyond}

When using more than three characters (e.g. 1, 2, \ldots, 6 with a die), the only formula which changes, in a rather obvious manner, is \eqref{Su19Apr1250}. For example, having three distinct possibilities for a symbol (say $\HH$, $\TT$, and $\RR$), we would now get
$$P_\cS = p_{{\sf H}}^{m_{\sf H}} \cdot p_{{\sf T}}^{m_{\sf T}} \cdot p_{{\sf R}}^{m_{\sf R}}$$
where $m_{{\sf H}}$, $m_{{\sf T}}$, and $m_{{\sf R}}$ is the number of symbols of each type found in $\cS$, and $p_{{\sf H}}$, $p_{{\sf T}}$, and $p_{{\sf R}}$ are the corresponding individual probabilities. The remaining formulas require no modification.


\begin{thebibliography}{9}
\bibitem{ref} G Blom and D Thornburn: {How many random digits are required
until given sequences are obtained?} \emph{J. Appl. Prob.} \textbf{19},
518-531 (1982)
\end{thebibliography}
\end{document}